\documentclass[12pt]{article}
\usepackage{amssymb}
\usepackage{latexsym}
\usepackage{exscale}

\begin{document}

\def\theequation{\thesection.\arabic{equation}}
\renewcommand{\rq}[1]{(\ref{#1})}
\newcommand{\ov}{\overline}
\renewcommand{\v}{ {\,\bf v }}
\newcommand{\w}{ {\,\bf w }}
\newcommand{\dist}{\mbox{dist}}
\newcommand{\dsize}{\displaystyle}
\newcommand{\supp}{\mbox{supp}}
\renewcommand{\Re}{\mbox{Re}}
\renewcommand{\Im}{\mbox{Im}}
\renewcommand{\supp}{\mbox{supp}}
\renewcommand{\div}{\mbox{div}}
\newcommand{\R}{{ \bf R  }}
\newcommand{\C}{ \bf C }
\newcommand{\plap}{\mbox{ $p$-Laplacian}}

\newcommand{\ra}{\rightarrow}
\newcommand{\BEL}{\begin{equation}\label}
\newcommand{\EE}{\end{equation}}
\renewcommand{\t}{\tilde}
\renewcommand{\medskip}{\vskip .5 cm}
\newtheorem{Thm}{Theorem}[section]
\newtheorem{Lemma}[Thm]{Lemma}
\newtheorem{Cor}[Thm]{Corollary}
\newtheorem{Prop}[Thm]{Proposition}

{

\title{Existence problems for the
$p$-Laplacian}
 \maketitle

\bigskip
\centerline{Julian Edward, Steve Hudson and Mark Leckband}

\bigskip

\noindent{\bf to be published in  Forum Mathematicum by De Gruyter}

\bigskip

{\centerline{\bf Abstract}} We consider a number of  boundary value
problems involving the $p$-Laplacian. The model case is $-\Delta_p
u=V|u|^{p-2}u$ for $u\in W_0^{1,2}(D)$ with $D$ a bounded domain in
${\bf R}^n$. We derive necessary conditions for the existence of
nontrivial solutions. These conditions usually involve a lower bound
for  a product of powers of the norm of $V$,  the measure of $D$,
and  a sharp Sobolev  constant. In most cases, these inequalities
are best possible. Applications to non-linear eigenvalue problems
 are also discussed.

\section*{ 1. Introduction}
\setcounter{section}{1} \setcounter{Thm}{0} \setcounter{equation}{0}

\vspace{5mm} Let $D$ be an open bounded region in
${\bf R}^n$, with $n\geq 1$. Define the $p$-Laplacian by
\begin{equation}\label{pLapdef}
\Delta_pu =  \div (|\nabla u|^{p-2}\nabla u) . \ \ \
\end{equation}
for $1<p<\infty$. Apart from its intrinsic interest, the $p$-Laplacian arises in the
study of non-Newtonian fluid mechanics both for $p\geq 2$ (dilatant
fluids) and $1<p<2$ (pseudoplastic fluids), see \cite{AS}. It also
arises in the study of
 of quasiconformal mappings \cite{Ho} and other topics in geometry \cite{U}.
 The one dimensional case also arises in models of turbulent flow of gas in
porous media \cite{OR}. In this work, we consider equations such as
\begin{equation}\label{pLap}
-\Delta_pu = V|u|^{p-2}u, \ u\in W_0^{1,p}(D),\ \ \
\end{equation}
where $V$ is assumed to be real-valued and integrable, and $u$ is
assumed to be complex-valued unless stated otherwise. We assume
nothing about the boundary of $D$. This paper generalizes previous
work by the authors \cite{DEHL}, (also see \cite{DH1},\cite{DH2}),
where the case $p=2$ is considered. Of course for $p=2$, $\Delta_p$
is the well-known linear Laplace operator.

Assuming a non-trivial solution to \rq{pLap}, we prove inequalities of the
form $$K\| V\| \geq 1.$$  To state a typical result, we need more notation.
Suppose $1<p<n$ and $q<\frac{np}{n-p}$.
Let $K = K_{q,p}=K(q,p,n,D)$ be the Sobolev constant of the embedding
$W^{1,p}_0(D)\to L^{q}(D)$. That is,
 \BEL{sob} K_{q,p}=\sup_{u\ne 0}  \frac{||u||_{q}}{ ||\nabla u ||_p}.\ \ \
 \EE
By Lemma \ref{SoboExtr}, equality is attained in \rq{sob} by a
nonnegative {\it extremal} function $u_*$. We have similar results
(see Lemma \ref{SoboExtr2}) and adopt similar notation when $p>n$
and $q\le \infty$. Our
first such result is:
\begin{Thm}\label{thmp0}  Suppose $1< p<n$ and
$\frac{n}{p}<r\leq \infty $. Let $q$ satisfy $1/r + p/q =1$, so that  $q<\frac{np}{n-p}$. Suppose
$V\in L^r(D)$, and $u\in W^{1,p}_0(D)$ is a nontrivial solution
of \rq{pLap}. Then
\begin{equation}
K^p_{q,p} ||V_+||_{r} \geq 1. \ \ \ \label{msrp0}
\end{equation}
 Let $u_*\ge 0$ be
an extremal for \rq{sob} given by Lemma \ref{SoboExtr}. Then
\rq{pLap} holds with $u=u_*$, and equality in \rq{msrp0} is attained when
$$ V(x)=  c_{*}u_*^{q-p}(x), {\ \  where\ \ } c_{*} = \| u_*\|_{q}^{-q}. \ \ \ $$
\end{Thm}

The proof of Theorem \ref{thmp0} is presented in Section 2, where we
also consider the case $p>n$ with $r\geq 1$, and the case $p<n$ with
$r\leq n/p$. We consider the case $p=n>1$ in Section 3. The
variations of Theorem \ref{thmp0} presented in the paper complement
results in \cite{OR} and \cite{BGG} and references therein, where
sufficient conditions for the   existence of solutions to \rq{pLap} are proved for $n=1$.

A key ingredient here is the imbedding of $W^{1,p}_0(D)$ into
various Banach spaces, based on inequalities of Sobolev or
Moser-Trudinger. An interesting theme here is that \rq{pLap} is the
Cauchy-Euler equation for the  extremals of these imbeddings. The
associated functions $V$ are extremals for our proclaimed lower
bounds, such as \rq{msrp0}. For a related discussion, see \cite{H}.
\medskip

We can estimate the Sobolev constant $K(q,p,n,D)$ in terms of $|D|$.
Let $B$ be a ball with $|B|=1$, and let $K^*_{q,p} = K(q,p,n,B)$.
Standard arguments involving scaling  and symmetrization, \cite{LL},
show that $K(q,p,n,D) \le |D|^{1/q -1/p +1/n}K^*_{q,p}$. Inserting
this into  \rq{msrp0}, we have \BEL{msrpD} |D|^{p/q -1 +p/n}
(K^*_{q,p})^p   ||V_+||_{(q/p)^*} \geq  1. \EE  This
is a minimal support result, which is a special type of unique
continuation result, see \cite{DEHL}. That is, if $|D|$ is too small
to satisfy \rq{msrpD}, then any solution to \rq{pLap} must vanish on
$D$. Similar remarks apply to many other results in this paper.

Equation \rq{pLap} can be viewed as a generalization of the non-linear eigenvalue equation: $ -\Delta_pu = E|u|^{p-2}u, \ u\in
W_0^{1,p}(D),\label{nonlinev}$ where $E\in \R$. Such eigenvalue
problems have been studied by a number of authors, see \cite{Le} and
references therein. A further generalization follows easily by
applying Theorem \ref{thmp0} to the function $V+E$:
\begin{Cor}\label{CorY1} Let
$p,q,r$ be as in Theorem \ref{thmp0}. Let $V\in L^r(D)$.
 Let  $u\in W^{1,p}_0(D)$ be a nontrivial solution of
\begin{equation}
 -\Delta_pu -V|u|^{p-2}u=E|u|^{p-2}u,\label{nonlinpot}
\end{equation}
 where $E\leq 0$ is a constant.   Then $ K_{q,p}^p
||(V+E)_+||_{r} \geq 1$.
\end{Cor}
\noindent This result  can by viewed as a lower bound on the
eigenvalues of \rq{nonlinpot}. Sufficient conditions for solvability of
\rq{nonlinpot} can be found in \cite{HR}. The case of Neumann
boundary conditions with $V\geq 0$, along with application to
degenerate parabolic equations, is studied in \cite{B}.

\medskip

The following corollary is basically a rephrasing of Theorem
\ref{thmp0}:
\begin{Cor}\label{CorY2} Let
$p,q,r$ be as in Theorem \ref{thmp0}. Let ${\cal G}=\{ V\in L^r(D):
V(x)\geq 0, \| V\|_r=1\}$.  Then among all pairs $(E ,V)\in {\bf
R}^+\times {\cal G}$ for which there exists non-trivial $u\in
W^{1,p}_0(D)$ such that
\begin{equation}
-\Delta_pu =EV|u|^{p-2}u, \label{nonlinev2}
\end{equation}     we have $$ E \geq \frac{1}{K^p_{pq} }.$$
Equality is attained when $u=u_*$, as defined in Lemma \ref{SoboExtr}, with $V  =u_{*}^{q-p}/\| u_{*}^{q-p}\|_r$.
\end{Cor}
 Note that ${1}/{K^p_{pq} }$ is the smallest eigenvalue of
\rq{nonlinev2}, and $u_*$ the corresponding eigenfunction.
 Corollary \ref{CorY2} can be compared with the following result in
\cite{CEP}. Fix $V_0\geq 0$ with $V_0\in L^{\infty}$ and
$||V_0||_r=1$. Let ${\cal  G}_1$ be the set of  all measurable
rearrangements
 of $V_0$. Then among all pairs $(E ,V)\in {\bf
R}^+\times {\cal G}_1$ for which there exists a positive $u\in
W^{1,p}_0(D)$ such that \rq{nonlinev2} holds,
 there exists $V_1$ that minimizes $E$.
  Furthermore, letting $u_1$ be
 the corresponding positive, normalized eigenfunction, there exists an increasing
 function $\phi$ such that $V_1=\phi (u_1)$.
A formula for $\phi$ is not given, and seems to be difficult to
deduce from the methods in \cite{CEP}.
 It would be
 interesting to see how $V_1,\phi$ relate to the pair
 $V,z\rightarrow z^{q-p}$ appearing in Corollary
 \ref{CorY2}.
 We remark  that analogues of Corollaries
 \ref{CorY1} and \ref{CorY2} hold for the variants of Theorem
 \ref{thmp0} that will follow.

Our paper is organized as follows. In Section 2, we present lower
bounds on $||V||$ based on the inequalities of Sobolev
and present examples to demonstrate sharpness. In Section 3, we
present lower bounds on an Orlicz norm of $V$, based on the
Moser-Trudinger inequality. In Section 4, we study two
generalizations of \rq{pLap}: the equation
$-\Delta_pu=V(x)|u|^{\beta}u$
 with $\beta \neq p-2$, and the equation
$-\Delta_pu=V(x)f(x,u,\nabla u)$. Finally, we have collected some
technical lemmas, perhaps not entirely new, in an appendix.

\section*{ 2. $L^r$ lower bounds.}
\setcounter{section}{2} \setcounter{Thm}{0} \setcounter{equation}{0}

In this section we assume that the potential $V$ belongs to some
Lebesgue space $L^r(D)$, with $r\ge 1$, and then show that $||V||_r$
must be large.

\subsection{The basic theorems}

We assume  $u\in W^{1,p}_0(D)$ and that \rq{pLap} holds in the
distribution sense. That is, for every $\psi \in C^{\infty}_{0}(D)$,
$$\int_D |\nabla u(x)|^{p-2}\nabla u(x) \cdot \nabla \psi(x)\;dx
=\int_D V(x)|u(x)|^{p-2}u(x)\psi(x)\;dx.$$ Define $\bar{q} = \frac{np}{n-p}$ for
$p<n$, and $\bar{q}=\infty$ for $p\geq n$.
 Recall that if $1 < p <
n$, then $u\in L^q(D)$ for $1\leq  q \leq \bar{q}$.  If $p> n$, then
$u\in C^{0}(D) \cap L^q(D)$ for $1\leq q \leq \infty$. Unless stated otherwise, we
assume $p\leq q $ and $r = (q/p)^*$ is the Holder conjugate of
$q/p$, meaning $1/r + p/q =1$, as in Thm. \ref{thmp0}.

In Theorem \ref{thmp}, we consider cases where \rq{msrp0} and
\rq{sob} have no extremal, for example, when $p<n$ and $q =
\bar{q}<\infty$. Note that in effect $r = (\bar{q}/p)^* = n/p$. In
this case, the Sobolev constant does not depend on $D$ and can be
computed explicitly \cite{T}. For brevity in proofs, we may
abbreviate it, letting $K = K_{q,p}$. For $p>n$, extremals do exist
for the Sobolev constant $K_{\infty ,p}$, but they don't solve
\rq{pLap} with $V\in L^1(D)$, see Section 2.2.
\medskip
\noindent
\begin{Thm} \label{thmp}
Let $1<p<n$ and let $u\in W^{1,p}_0(D)$ be a nontrivial solution of
\rq{pLap}, with $V\in L^{n/p}(D)$. Then\begin{equation}
K_{\bar{q},p}^p ||V_+||_{n/p}
> 1.\ \ \ \label{msrp}
\end{equation}
If $p>n$
 and $V\in L^1(D)$, then
\begin{equation}
K_{\infty,p}^p ||V_+||_1> 1.\label{msrp3}
\end{equation} The constant
1 is sharp in both \rq{msrp} and \rq{msrp3} when $D$ is a ball.
\end{Thm}
If desired, one can replace $||V_+||_{n/p}$ in \rq{msrp} or \rq{msrp3} by $||V||_{n/p}$.
In that case the constant 1 is still sharp, with the same proof.

A maximal principle formulated in (\cite{GT}, Theorem 10.10) implies
an inequality of the form $c\| V_+\|_{n/p}>1$, similar to \rq{msrp}.
But it is not clear whether the constant $c$ is sharp.

The following extension of Theorem \ref{thmp0} handles the cases
where $V\in L^r$ with $r>n/p$ for $p<n$,
 and $r>1$ for $p>n$. In both cases, the Sobolev embedding is compact,
 and we have extremals for \rq{sob} satisfying \rq{pLap} with $V$ in a
Lebesque space.
\medskip

\begin{Thm}\label{thmp2}  Suppose
$\max (1,\frac{n}{p})<r\leq \infty $. Suppose $V\in L^r(D)$,
and $u\in W^{1,p}_0(D)$ is a nontrivial solution of \rq{pLap}.
 Then
\begin{equation}
K^p_{q,p} ||V_+||_{r} \geq 1. \ \ \ \label{msrp2}
\end{equation}
 Let $u_*\ge 0$ be
an extremal for
\rq{sob} given by Lemmas \ref{SoboExtr2}
 or \ref{SoboExtr}. Then \rq{pLap} holds, and equality in \rq{msrp2} is attained, when
\BEL{V} V(x)= c_{*}u^{q-p}_*(x), {\ \  where\ \ } c_{*} = \frac{
1}{\| u_*\|_{q}^{q}}. \ \ \ \EE
\end{Thm}
\medskip

We have the following simple consequence of Theorems ~\ref{thmp} and
~\ref{thmp2}, which includes Corollary \ref{CorY1} as a special
case.
\begin{Cor}\label{CorY}  Let  $u\in W^{1,p}_0(D)$ be a nontrivial solution of
$$-\Delta_pu -V|u|^{p-2}u=E|u|^{p-2}u.$$ where $E\leq 0$ is a constant.
\medskip
 (a) If $p<n$ and $V\in L^{n/p}(D)$, then $K_{\bar{q},p}^p ||(V+E)_+||_{n/p} > 1.$

 (b) If $p>n$ and $V\in L^1(D)$, then $K_{\infty,p}^p ||(V+E)_+||_1 > 1$.

 (c) If $\max (1,\frac{n}{p})< r\leq \infty $ and $V\in
L^r(D)$, then $K_{q,p}^p ||(V+E)_+||_{r} \geq 1$.
\end{Cor}
 The obvious analogue of Corollary \ref{CorY2} also holds.
\medskip

\noindent {\it Proof of Thm \ref{thmp2}:} Since $\max
(1,\frac{n}{p})<r\leq \infty $ we have $q\leq \bar q$, which allows
Sobolev's inequality \rq{sob2} below. Also, using Green's identity
(see Lemma ~\ref{GreenG}), and Holder's inequality based on $1/r +
p/q =1$,
\begin{eqnarray}
\| u\|_{q}^p & \leq & K^p \int_{D}| \nabla u|^{p}dx\label{sob2}\\
& = & K^p \int_{D}V |u|^p dx \nonumber \\
& \leq & K^p \int_{D}V_+|u|^{p}dx\label{posi} \\
& \leq  & K^p \| V_+\|_{r}\| u\|_{q}^p. \label{holder2}
\end{eqnarray}
Since $\| u\|_q>0$, we have $K^p \| V_+\|_{r}\geq 1$, and \rq{msrp2}
holds.  This type of proof, which also appears in \cite{DEHL} and
later in this paper, will be referred to informally as a {\it
minimal support sequence}.
\medskip
We now show that equality can be attained in \rq{msrp2}.
By Lemmas \ref{SoboExtr2}, \ref{SoboExtr} and \ref{EuLap},
there is a $u_*\geq 0$ for which \rq{sob2} holds with equality, with $\| \nabla u_*\|_p=1$ and
 \BEL{PDE2} -\Delta_p u_*=
\frac{u_*^{q-p}}{ \| u_*\|_{q}^{q}}u_*.
 \EE
Setting  $V=u_*^{q-p}/\|u_*\|_q^q$, we also have equality in
\rq{posi}. Equality occurs in \rq{holder2} because
$(|u_*|^p)^{\frac{q}{p}-1}=|u_*|^{q-p}=cV_+$, see for example
\cite{LL}, p.45. $\Box$ \noindent

For more insight into why the extremal for Sobolev's inequality
generates the extremal for \rq{msrp2}, the reader is referred to
\cite{H}.
\medskip

\noindent {\bf Remark:} the method of proof of Theorem \ref{thmp2}
can easily be applied to study the equation $\div (a(x)|\nabla
u|^{p-2}\nabla u)=V|u|^{p-2}u$, where $a$ is a positive function
with $a,1/a\in L^{\infty}$. Other results in
this paper can also be generalized in this way.

\medskip

\noindent {\it Proof of Thm \ref{thmp}:} We begin with the case
$p<n$. The proof of Thm \ref{thmp2} shows that $K^p \| V_+\|_{r}\geq
1$ when $r = (q/p)^* = (\bar{q}/p)^* = n/p$. It is well-known
\cite{T} that when $q=\bar q$, and $u$ is nontrivial, equality
cannot occur in \rq{sob2} for any bounded domain $D$.   So, $K^p \|
V_+\|_{n/p}> 1$, as desired.

To show that the constant 1 is sharp in \rq{msrp}, we will set
$D=B_{\hat{R}}(0)$ and construct $u,V$ on $D$ so that
$K^p||V||_{n/p} \to 1$ as $\hat{R}\to \infty$. Since $1 <
K^p||V_+||_{n/p} \le K^p||V||_{n/p}$, this
 implies $K^p||V_+||_{n/p}$ can be arbitrarily close to 1.
Let $\rho = |x|$. For radial $u$,
\begin{equation}
\Delta_p u(\rho) = |u_{\rho} |^{p-2} ( (p-1) u_{\rho \rho} +
\frac{n-1}{\rho} u_{\rho})= (p-1)|u_{\rho} |^{p-2} ( u_{\rho \rho} +
\frac{s-1}{\rho} u_{\rho}),\label{polar}
\end{equation}
for $s=\frac{n-1}{p-1} + 1$.  The extremal for the Sobolev embedding
$W^{1,p}_0({\bf R}^n)\ra L^q({\bf R}^n)$ for $1<p<n$ with critical
index $q=pn/(n-p)$ is given by
\begin{equation}
v(x)=(1+\rho^{p'})^{(p-n)/p},\label{talentiext}
\end{equation}
where
$p'=p/(p-1)$, see \cite{T}. Define $V_{v}$ by $-\Delta_{p}v = V_{v}
v^{p-1}$. Applying the minimal support sequence to $v,V_v$ with
$r=n/p$, and with $D={\bf R}^n$ (temporarily), we get
$K^p||V_v||_r=1$.

Below, let $C$ be a positive constant whose value may change at each
step, and let $R$ be a constant that eventually will approach
infinity, so we can assume without loss of generality that $R>C$.
The other constants below may depend on $R$, but $K$ does not, as it
is independent of dilation in this critical case. Let $\hat{R}
> R+1$ to be specified later. Set
$$u=\cases{
v,& $0 \leq \rho < R$,\cr
 a - b \rho,& $R \leq \rho < R+1$, \cr
c\rho^{2-s} + d,& $R+1 \leq \rho \leq \hat{R}$.\cr}$$
To make $u$ and $u_{\rho}$ continuous at $\rho = R$ and at
$\rho = R+1$, let
\begin{eqnarray*}
a - bR &=& v(R),\\
 b &= & \frac{(n-p)}{(p-1)}R^{p'-1} (1+R^{p'})^{-n/p},\\
 c &= &  (R+1)^{s-1} R^{p'-1}(1+R^{p'})^{-n/p},\\
 d & =&  - \frac{(n-p)}{(p-1)}R^{p'-1}(1+ R^{p'})^{-n/p}.
\end{eqnarray*}

 Note that by \rq{talentiext},
$b=|v_{\rho}(R)|\leq Cv(R)/R$ for some $C$ independent of $R$.
Assuming $R>C$,
\begin{equation}
u(R+1) = a-b(R+1) = v(R)- b\ge v(R)(1-C/R) >0
\label{error}\end{equation} Since $\lim_{\rho\to \infty} c\rho^{2-s}
+ d =d <0$ and $u$ is continuous, there is some $\hat{R} > R+1$ so
that $u(\hat{R}) = 0 $. Set $D=B_{\hat{R}} (0)$, and note that $u
\subset W^{1,p}_0 (D)$ is radial and nonnegative.
 Define $V_u $ by $ -\Delta_p u = V_u u^{p-1}$. For $R \leq \rho < R+1$, we get $|V_u| = |\frac{(n-1)
b^{p-1}}{\rho (a-b\rho )^{p-1}}| \leq C R^{-p}$. For $R+1 \leq
\rho$, we have
\begin{equation}
\Delta_p(\rho^{2-s})=0,\label{harm}
\end{equation}
hence
 $V_u =0$. So,
$$ \int _{D} |V_u|^{n/p} \le C \int _{R \leq \rho \leq R+1} R^{-n} + \int _{B_{R} (0)} |V_v|^{n/p} .$$
As $R\rightarrow \infty$, the first integral is bounded
by $CR^{-1}$ which converges to $0$. The second integral converges to
$||V_v||_{n/p}^{n/p} = K^{-n}$, so $K^{p}\|V_u \|_{L^{n/p}(D)}
\rightarrow 1$, as desired. This proves sharpness of \rq{pLap} when
$p<n$.

\medskip
We now address the case $p>n$. The proof of Thm \ref{thmp2} applies,
with $q = \bar q = \infty$ and $r=1$, and hence $K^p\|V_+\|_{1} \geq
1$. We will show in the next subsection, in Thm \ref{Steve1}, that
the last inequality is strict, ie. $K^p\|V_+\|_{1} > 1$, and that
the constant 1 cannot be improved, completing the proof.  $\Box$

\subsection{Inequalities for bounded solutions with $V\in L^1$}

\medskip
If $p>n$ then $u\in L^{\infty}$ and we may consider $V\in L^1(D)$ in
Theorem \ref{thmp}.  The lower bound for $||V_+||_1$ still holds,
but the usual Euler-Lagrange equation \rq{PDE} does not, raising
interesting new questions about sharpness and extremals. We prove an
analogue of Theorem \ref{thmp} replacing $L^1(D)$ with the space $M$
of signed measures $V$ on $D$, (see e.g. \cite{Ru} for the
definition and properties of signed measures) with norm $\|
V\|_M=|V|(D) < \infty $. In the special case where $V\in L^1(D)$, we
have $\| V\|_M=|V|(D)=\int_D|V(x)|dx=\| V\|_1.$ We also recall
$V_+=\frac12(V+|V|)$. Assume \rq{pLap} holds for $V\in M$, ie.
$$\int_D|\nabla u|^{p-2}\nabla u \cdot \nabla \phi =<V,|u|^{p-2}u\phi >, \ \forall \phi \in C_0^{\infty}(D).$$
Here $<*,*>$ denotes the pairing of signed measures with continuous
functions. Let $\phi_n \in C_0^{\infty}(D)$ converge in
$W^{1,p}(D)$, and hence in $C^0(D)$, to $ \ov u$. It then follows
that
\begin{equation}
\int_D|\nabla u|^{p}  =<V,|u|^{p} >.\label{Greenm}
\end{equation}
\bigskip
\begin{Thm}\label{Steve1}
Assume $u\in W_0^{1,p}(D)$ is a nontrivial solution of \rq{pLap},
with $V\in M$. Then \BEL{rBe1}K_{\infty,p}^{p}||V_+||_M \ge 1. \EE
Equality is attained when $V=K^p_{\infty,p}\ \delta_z(x)$, where
$\delta_z$ is a Dirac mass at some point $z\in D$. When $V\in
L^1(D)$, equality is not possible, but the constant 1
cannot be replaced by any larger constant.
 \end{Thm}
\noindent Proving the last assertion of this theorem will complete
the proof of Theorem ~\ref{thmp}.

\medskip

 \noindent {\it Proof:} Let $K=K_{\infty,p}$. Using
\rq{Greenm}, we have the minimal support sequence
\begin{eqnarray}
\| u\|_{\infty}^p& \leq & K^p\| \nabla u\|_p^p\nonumber \\
& =& K^p<V,|u|^p>\nonumber \\
& \leq & K^p<V_+,|u|^p>\nonumber \\
& \leq & K^p \|V_+\|_M\|u^p\|_{\infty}.\label{msrmeasure}
\end{eqnarray}
The inequality \rq{rBe1} follows immediately. We now demonstrate extremals for this inequality.
By Lemma ~\ref{SoboExtr2}, there is a non-negative extremal $u_* \in W^{1,p}_0(D)$ for
the Sobolev inequality used in \rq{msrmeasure}, that is
\begin{equation}
||u_*||_\infty = K||\nabla u_*||_{p}\ . \label{kinfty1}
\end{equation}
This $u_*$ represents a scalar multiple of the extremal in Lemma ~\ref{SoboExtr2}. That is, we do not assume  $||\nabla u_*||_{p}=1$.
Since $u_*$ is continuous, it attains its maximum value at some
point $z\in D$. So, \rq{kinfty1} can be rewritten
$$
K = \frac{|u_*(z)|}{||\nabla u_*||_p}.
$$
where $u_*$ maximizes the right-hand side among all $u\in
W^{1,p}_0(D)$.  We now apply the Euler-Lagrange method to this. Let
$\phi \in C_0^{\infty}(D)$, and let
$u_{\epsilon}(x)=u_*(x)+\epsilon\phi(x)$. Then
$$0=\frac{d}{d\epsilon}(\frac{|u_{\epsilon}(z)|}{||\nabla
u_{\epsilon}||_p})|_{\epsilon=0} $$ gives $$0 =\phi (z)-u_*(z)\|
\nabla u_*\|_p^{-p}\int_D|\nabla u_*|^{p-2}\nabla u_* \cdot \nabla
\phi dx.$$ Normalizing, we can assume $u_*(z)=1$, so
$1=K\|\nabla u_*\|_p$, and the above shows $u_*$ is a weak solution
to \BEL{weak} -\Delta_pu=K^{-p}\delta_z, \EE and it satisfies
\rq{pLap} with $V_*=K^{-p}\delta_z$. Since $||\delta_z||_M = 1$, we
have $K^p||V_*||_M =1$, and hence $V_*$ is an extremal.
\medskip

Now suppose \rq{rBe1} was an equality for a $V\in L^1(D)$ with
corresponding solution $u$ to \rq{pLap}. Then the argument leading
to \rq{msrmeasure} would imply $\| u\|_{\infty}^p = K^p\| \nabla
u\|_p^p$, so $u$ would be an extremal for the Sobolev inequality.
Hence, by the argument above, $-\Delta_pu=K^{-p}\delta_z|u|^{p-2}u = V |u|^{p-2}u$ as distributions, for some $z\in D$.
This implies $V |u|^{p-2}u = 0$ a.e., so that $\Delta_pu =0$ as a distribution.
But then $u\equiv 0$, a contradiction.

\medskip
Next, to show that the constant 1 in \rq{rBe1} is sharp for $V\in L^1$,
we construct examples on $D=B_1(0)$ such that $K^{p}||V||_1 \to 1$.
Define $u_*$ as a Sobolev extremal as above. By standard
symmetrization, we may assume $u_*$ is radially non-increasing, so
that $z=0$, and by \rq{weak} $-\Delta_pu_* \equiv 0$ away from $0$.
Using \rq{polar} and solving an ordinary differential equation, we
get $u_*(x)= 1 -\rho^{\frac{p-n}{p-1}}$.

We approximate $u_*$ by a function $u$ such that
 $V_u =\frac{\Delta_pu}{u^{p-1}} \in L^1$. Let $\epsilon >0$ and let
$u=u_*$ outside $B_\epsilon(0)$. On $B_\epsilon(0)$, define $u(\rho
) = a - b\rho ^2$ such that $u'(\epsilon)$ exists. Note $\Delta_p u
\in L^1(B_1(0))$, so $V_u \in L^1$ too. Also, $\Delta_p u \le 0$
there, because $p>n$, $b>0$, and $u'' + (n-1)u'/((p-1)\rho) = -2b[1
- (n-1)/(p-1)] <0$. So $V_u\ge 0$, and
$$||V_u||_1 = \int_{D} V_u \ dx \le u_*^{1-p}(\epsilon) \int_{D} V_u u^{p-1} \ dx
= -u_*^{1-p}(\epsilon) \int_{D} \Delta_p u \ dx.
$$ Let $v\in
C^\infty_0(D)$ with $v\equiv 1$ on $B_\epsilon(0)$. Then
$$\int_{D} \ \Delta_p u \ dx = \int_{D} v\ \Delta_p u \ dx =
\int_{D} |\nabla u|^{p-2} \nabla u\ \cdot \nabla v\  dx = \int_{D}
|\nabla u_*|^{p-2} \nabla u_*\ \cdot \nabla v\  dx. $$ By \rq{weak}
and the definition of $\Delta_p u_*$, we get $\int_{D} \ \Delta_p u
\ dx = -K^{-p}$. As $\epsilon \to 0$, we have $u_*(\epsilon)\to 1$,
so that $||V_u||_1 \to K^{-p}$, as desired.

 \subsection{On the failure of other norms on V}

The proofs of Thms. \ref{thmp} and \ref{thmp2} depend on the Sobolev
and Holder inequalities, which impose the restriction $r\ge n/p$
(when $p<n$). In this section, we prove that if $ r < n/p$,
$||V||_r$ can be arbitrarily small.

\begin{Thm}\label{NewCn}   Let $p<n$ and $D=B_1(0)\subset \R^n$. For every $1\le r<\frac np$,
 and every   $\delta >0$, there is a potential $V_{\delta}$, with a nontrivial solution
  $u \in  W^{1,p}_0 (D) $ of \rq{pLap}, with $\|V_\delta \|_r < \delta $.
\end{Thm}
\noindent {\it Proof.}  We will specify $\epsilon =\epsilon (\delta)\in ( 0,1/2) $  later. With the usual convention that
$u(x) = u(|x|) = u(\rho)$, define
$$ u(\rho)  =\cases{ a-b\rho^{\frac{p}{p-1}}, & $ 0 \leq \rho < \epsilon.$ \cr
\rho ^{2-s} - 1, & $ \epsilon \leq \rho \leq 1,$ \cr}$$
where $ s = \frac{n-1}{p-1} + 1$, and $a,b$ are chosen to make $u\in
C^1_0(D)$. So, $b =\frac{n-p}{p}\epsilon^{-\frac{n}{p-1}}$ and $a=
\frac{n}{p} \epsilon^{2-s} -1$. Since $u$ is radial, by \rq{polar}
we have
$$ \Delta_p u(\rho) =  (p-1)|u_{\rho} |^{p-2} ( u_{\rho \rho}
+ \frac{s-1}{\rho} u_{\rho}).$$

 On $\epsilon \leq \rho \leq 1$, by \rq{harm} we have
$ \Delta_p u(\rho) = 0$ and $V_{\delta}\equiv 0$. Let $C$ denote
 positive constants that vary from line to line. For $0 \leq \rho <
\epsilon$, $ u_{\rho}(\rho) = -\frac{p}{p-1}b\rho ^{\frac{1}{p-1}}$,
$ u_{\rho \rho}(\rho) = -\frac{p}{(p-1)^2}b\rho^{\frac{1}{p-1}-1}$
and $-\Delta_p u(\rho)=C \epsilon^{-n}$; hence $|V_{\delta}| =
\frac{C\epsilon^{-n}}{|u|^{p-1}} \leq C \epsilon^{-p}$. So,
$\|V_{\delta} \|_{r}^{r} \leq C\epsilon^{n-rp} < \delta$, for small
enough $\epsilon$.

\bigskip

\section*{ 3. The case $p=n$: Orlicz lower bounds.}
\setcounter{section}{3} \setcounter{Thm}{0}
\setcounter{equation}{0}\setcounter{subsection}{0}

\subsection{The critical case $p=n\geq 2$ : $V \in L \log^{n-1} L(D)$}

When $p=n$, the proof of Theorem \ref{thmp2} does not extend to $r =
\frac np =1$, because $W^{1,n}_0(D)$ does not embed into
$L^{\infty}(D)$. In this section we assume instead that $V$ is in
the Orlicz space $L\log^{n-1} L(D)$, so that $\int_D
|V|(\log^{n-1}(1+ |V|))dx $ is finite, and prove an analogue of
Theorem ~\ref{thmp2} for that space. By Theorem \ref{LlogkL} no such
analogue holds for $V\in L\log^k L(D)$, with $0\le k < n-1$, and
therefore not for $V\in L^1$.

As a substitute for the Sobolev inequality we will use the Moser-Trudinger
inequality (see \cite{M}). Let $\alpha_n=
(n^{n-1}\omega_{n})^{1/n}$, where $\omega_{n}$ is the surface area
of the unit sphere in $\bf{R^n}$. Suppose $u\in W^{1,n}_0(D)$ is real-valued.
 The inequality is
$$\int_D\exp ((\alpha_n |u(x)|/\|
\nabla u\|_n)^{\frac{n}{n-1}})\ dx\leq C_n|D|.$$
Let $0< \alpha< (\alpha_n)^{n}$ be a fixed
constant. Define
\begin{equation}
M(t) =  \int_{0}^{\alpha t}e^{s^{\frac{1}{n-1}}} - 1\; \
ds\label{Mn}
\end{equation}
and
$$N(s) =  \int_{0}^{s/\alpha} \log^{n-1}(t+1)\;dt\ .$$
These are complementary Orlicz functions, see \cite{KR}. This
non-standard choice for $M,\ N$  allows
  an explicit formula for $N(s)$, which is useful because $N$ is used
in the definition of the norm of $V$.
Let $P_0(x)=1$, and $P_m(x) =\sum_{k=0}^{m} \frac{(-1)^k
m!}{(m-k)!}x^{m-k}$ for $m \geq 1$. Then, for $n \geq 2$,
 $$N(s)=(1+\frac{s}{\alpha} )P_{n-1}(\log(1+\frac{s}{\alpha}))+(-1)^{n}(n-1)! \ .$$
The Orlicz class $L^{N}(D) $ is the set of measurable $u$ such that
$\int_D N(|u|)dx <\infty$ with $L^{M}(D)$ defined similarly. Fix
$c$, with $\alpha^{1/(n-1)}<c<(\alpha_n)^{n/(n-1)}$. By \rq{Mn}
there is a $C>0$ such that both $M(t),\;M'(t)<Ce^{
ct^{\frac{1}{n-1}}}$ for all $t>0$. So the Moser-Trudinger
inequality gives
\begin{equation}
     \int_D  M\left ( \frac{| u(x)|^n}{\|\nabla u\|^n_n}\right )   dx \leq K_{M}|D|,\ \
      \forall u \in W_{0}^{1,n}(D),  \label{MTII}  \end{equation}
where the optimal constant $K_{M}$ depends on $D$, but is dilation-invariant and independent of $u$.

\noindent {\bf Remark:} An example in \cite{M} shows the \rq{MTII}
does not hold if $\alpha= (\alpha_n)^{n}$ for $n \geq 3$. The case
$n=2$ is discussed later in a remark below.

\bigskip

 We define the norm of $V\in L\log^{n-1} L(D)$ by
 \BEL{lux}\|V\|_{N} = \inf \left\{ \lambda
+\frac{\lambda}{K_{M}|D|} \int_D
N\left(\frac{|V(x)|}{\lambda}\right) dx;\ \lambda
>0 \right\} <\infty , \label{defN}\EE
see \cite{KR}. For fixed $V$, we set $F(\lambda) = \lambda \int_D
N\left(\frac{V_+(x)}{\lambda}\right) dx$, so that $\|V_+\|_{N} =
\inf \left\{ \lambda + \frac{F(\lambda)}{K_{M}|D|}\right\}$. With
the norm for $L^M$ defined analogously, standard arguments show that
the injection $W^{1,n}_{0}(D)\rightarrow L^{M}(D)$ is compact.

\begin{Lemma}\label{EuLan2}
There exists a non-negative extremal $u$ for \rq{MTII}. Furthermore,
the Euler Lagrange equation for \rq{MTII}, with  the normalization
$||\nabla u||_n=1$, is $$ -\Delta_n u =V u^{n-1} $$ where \BEL{V2}
V(x)=  \frac{M'(u^{n}(x))}{\int_D M'(u^{n}) u ^{n} dx}\in
L\log^{n-1}L(D). \EE
\end{Lemma}

\noindent {\it Sketch of Proof:} By the compactness of
$W^{1,n}_{0}(D)\rightarrow L^{M}(D)$, there is an extremal for
\rq{MTII}.  Variational work similar to that done for Theorem
\ref{Steve1} gives the equation in \rq{V2}. That $V\in
L\log^{n-1}L(D)$ follows from \rq{Mn}.
\bigskip

Our main result for the case $p=n$ is:
\begin{Thm}\label{Orl}Assume that \rq{pLap}
has a nontrivial solution $u$ for $V\in L\log^{n-1} L(D)$. Then
\BEL{ThmOrl} \dsize K_{M}|D| \| V_+\|_{N}\geq 1, \EE where $K_{M}$
is the optimal constant in \rq{MTII}. With $u=u_*$ and $V$ as in
Lemma \ref{EuLan2}, equality is attained in \rq{ThmOrl}, and
\rq{pLap} holds.
 \end{Thm}

\noindent {\bf Remark:} The method of proof of Theorem \ref{Orl} can
be adapted to the Orlicz function $$\tilde{M}(t) =e^{(\alpha_n^n
t)^{1/(n-1)}} -  \sum_{k=0}^{n-1} (\alpha_n^n t)^{k/(n-1)}/k! \
,$$for which, by \cite{Li}, \rq{MTII} has extremals for $n\geq 2$.
For $n=2$, we observe this formula is (\ref{Mn}) with $\alpha =
(\alpha_2)^{2}$.

\bigskip

Theorem \ref{Orl} follows immediately from the following:
\begin{Thm}\label{Orl2}
Suppose that  \rq{pLap} has a nontrivial solution $u$ with $V\in
L\log^{n-1} L(D)$. Then, for every $\lambda>0$,

\BEL{inf}  \lambda K_{M} |D| + F(\lambda) \geq 1.
  \EE
Equality can be attained in \rq{inf} with $u_*$ and $V$ as in Lemma \ref{EuLan2}, and \rq{pLap} holds.
\end{Thm}

\medskip
\noindent {\it Proof of Theorem \ref{Orl2}.} Fix $u,V$ with $\|
\nabla u\|_n =1$. Let $U =|u(x)|^{n}$. For fixed $\lambda >0$, set
$v = \frac{V_+(x)} \lambda$.  It is well known for any Orlicz pair
$(M,N)$ that Young's inequality gives \BEL{Young}Uv \leq M(U) + N(v)
\EE with equality if and only if $v = M'(U)$. By  Green's identity
(Lemma \ref{GreenG}), the definition of $U$ and \rq{Young},
\begin{eqnarray*}
1=\|\nabla u \|^n_n &=&
 \int_D |u|^n V dx\\
&\leq &\int_D U V_+ dx\nonumber \\
& \leq &   \lambda \int_D M(U)dx + F(\lambda) . \label{minF}
\end{eqnarray*}
By \rq{MTII},  $\dsize \int_D M(U)dx  \le K_{M}|D|$ and \rq{inf}
follows.

Let $u_*$ and $V$ be as in Lemma \ref{EuLan2}, so that
 $-\Delta_{n} u_*= Vu_*^{n-1}$, where $V = \omega^{-1}M'(u_*^n) \ge
0$, where
$${\omega} = { \int_D M'(u_*^n)u_*^n   dx }.$$
Let $U = U_* = u_*^{n}$, so $\int_D M(U_*)dx = K_{M}|D|$ (see \rq{MTII}).
Setting $\lambda = \omega^{-1}$, we have $v=V_+/\lambda=M'(U_*)$,
so equality holds in \rq{Young}. From the definitions of $V$, $U_*$ and $\omega$, we have $\int_D|U_*V|dx = 1$.
Then, integrating \rq{Young},
$$1 = \lambda \int_D M(U_*)\ dx + F(\lambda) = \lambda K_{M}|D| +  F(\lambda).$$ Thus for
these choices of $u$, $V$ and $\lambda $,  \rq{inf} is an equality,
and also \rq{pLap} holds. $\Box$.
\medskip

\begin{Cor}\label{CorO}
Let $u\in W^{1,n}_0(D)$, and $V\in L\log^{n-1}L(D)$. If
$$-\Delta_nu -V|u|^{n-2}u=E|u|^{n-2}u,$$
with $E\leq 0$, then
$$ K_{M}|D|\ \| (V+E)_+\|_{N}\geq
1.$$
\end{Cor}

\begin{Cor}\label{CorY3}
Let $M,N$ be as above. Let ${\cal G}=\{ V\in L_N(D): V(x)\geq 0, \|
V\|_N=1\}$.  Then among all pairs $(E ,V)\in {\bf R}^+\times {\cal
G}$ for which there exists non-trivial $u\in W^{1,n}_0(D)$ such that
$$
-\Delta_nu =EV|u|^{n-2}u,
$$
    we have $$ E \geq
\frac{1}{K_{M}|D|},$$ with equality attained by $u$ and $V/E$, with
$u,\;V$ as in Lemma \ref{EuLan2}. Furthermore, $V(x)=\phi (u(x))$,
with $\phi$ an increasing function whose explicit formula can be
found in Lemma~\ref{EuLan2}.

\end{Cor}
\noindent As with Corollary \ref{CorY2}, this result can  be
compared  to the result in \cite{CEP} cited in the introduction.

\subsection{A counterexample for $V\in L\log^{k}L(D)
,\; k < n-1$.}

The purpose of this subsection is to present:
 \begin{Thm}\label{LlogkL}  Let $n\geq 2$ and $0\le k<n-1$.  Let $N(s)=\int_{0}^{s/\alpha} \log^{k}(t+1)\;dt$. For every $\delta> 0$, we can find a non-negative $V_{\delta}\in L\log^{k}L(B_1(0))$, and a positive solution
 $u\in W^{1,n}_{0}(B_1(0))$ of $-\Delta_n u = V_{\delta} u^{n-1}$, such that
  $ \|V_{\delta}\|_{N}<\delta $.
\end{Thm}

\medskip
\noindent Note that when $k=0$,  $V_\delta \in L^1$.

\medskip

 \noindent
 {\it Proof of Theorem ~\ref{LlogkL}:} Our constructed functions will be radial and positive. By
\rq{polar},
$$ \Delta_n u(\rho ) =  (n-1)|u_{\rho}|^{n-2} ( u_{\rho \rho} + \frac{1}{\rho} u_{\rho}).$$
Let $\delta >0$ be given and $ 0< \epsilon < 1/2 $ to be determined later.
Let  $$ u(\rho)=\cases{a-b\rho^{\frac{n}{n-1}}
   & if $0 \leq \rho < \epsilon $ ,\cr
  -\log (\rho)  & if $\epsilon \leq  \rho \leq 1$  \cr}$$
where $a$ and $b$ are chosen below so that $u$ is differentiable.
Note that $\Delta_n u(\rho)=0$ for $\epsilon \leq  \rho \leq 1$.
Continuity at $ \rho = \epsilon$ of $u_{\rho}$ requires $ b =
\frac{n-1}{n} \epsilon^{-\frac{n}{n-1}}$, and of $u$ requires $a =
\frac{n-1}{n} - \log(\epsilon )$. Let $C$ denote a constant which
may change from line to line. We define $V_{\delta}=V$ by the
equation $-\Delta_n u = V u^{n-1} $, which gives $ V=
   0$ for  $\epsilon \leq  \rho \leq 1$, and $V
   \leq \frac{C\epsilon^{-n}}{u^{n-1}(\epsilon)}$  for $0 \leq \rho <
   \epsilon$. Hence,
$$\int_{B_1(0)} N(|V|)\;dx \leq \frac{C|\log (\epsilon)|^{k}}{|\log (\epsilon)|^{n-1}}.$$
Let $\lambda = \delta/2$.  Observe that $\log(at+1)\leq a\log(t+1)$
for $a>1$. So, by the integral definition of $N(s/\lambda)$, we have
$N(s/\lambda)\leq \lambda^{-(k+1)}N(s)= CN(s)$, for all $s>0$. Hence
$$\lambda + \frac{\lambda}{K_M|{B_1(0)}|}\int_{B_1(0)} N(|V|/\lambda)\;dx\leq \delta/2+ C\int_{B_1(0)} N(|V|)\;dx \leq \delta/2+\frac{C|\log (\epsilon)|^{k}}{|\log (\epsilon)|^{n-1}}.$$
Choosing $\epsilon$ so that $C|\log (\epsilon)|^{k-(n-1)}<\delta/2$,
the result follows from \rq{defN}.$\Box$

\section*{ 4. Equations with other nonlinear terms}
\setcounter{section}{4} \setcounter{Thm}{0} \setcounter{equation}{0}

\medskip
We consider the equation \BEL{equbeta1} -\Delta_p u = V|u|^{\beta} u
\EE where $\beta \geq -1$. This is assumed in the
weak sense, that
$$\int_D  |\nabla u|^{p-2} \nabla u \cdot \nabla \psi \,dx= \int_D V|u|^{\beta} u\psi\,dx ,\
\forall \psi \in C_0^{\infty}(D).$$

\medskip
\begin{Thm}\label{Sbeta} Let $1<p<\infty$ and $r>1$. Define $\hat {q} = r(\beta+2)/(r-1)$, and
assume $\hat{q} \leq \ov q $. Let $u\in W^{1,p}_0(D)$ be a nontrivial weak
solution of \rq{equbeta1} with $V\in L^r(D)$. Then \BEL{betaS} K^p
\|V_+ \|_{r}\ ||u||_{\hat {q}}^{\beta +2-p} \geq 1. \EE where
$K=K_{\hat{q}p}(D)$. If $\hat q < \ov q$, equality can be attained
in \rq{betaS}.
\end{Thm}

\noindent {\it Proof:} We prove the result for $p\neq n$; the proof
for $p = n$ is almost identical.
 By Sobolev's inequality, $u\in L^{\ov
q}(D)$, and hence $u|u|^{\beta}V\in L^{(\bar{q})^*}(D)$. It follows by Lemma~\ref{GreenGG} and Holder's
inequality that
\begin{eqnarray} \label{nonlinS} ||u||_{\hat q}^p \leq K^{p} ||\nabla u||_p^p
\leq K^{p}\int_D |u(x)|^{\beta +2}V_+(x) \, dx
 \le K^p||u||_{\hat q}^{\beta +2}||V_+||_{r}
\end{eqnarray}
from which \rq{betaS} follows.

If $\hat q < \ov q$, Lemma
~\ref{EuLap}  provides a $u_* \ge 0$, with
 $||\nabla u_*||_p =1$,
such that  $ -\Delta_p u_* =  c u_* ^{\hat q -1}$ with
$$c = \|u_*\|_{\hat q}^{-\hat q} = K^{-p}\|u_*\|_{\hat q}^{p-\hat q}.  $$
So, $-\Delta_p u_* = Vu_*^{\beta +1}$, which is \rq{equbeta1}, with
$V =V_+= cu_*^{\hat q -2 -\beta}$. Thus, $\|V_+\|_{r} = c\| u_*
\|_{\hat q}^{\hat q - \beta-2}$, which gives equality in
\rq{betaS}. $\Box$

\medskip

We now consider equations such as $-\Delta_pu=V(x)|\nabla u|^{p-1}$, and give conditions under which  $K_{p,\ov{q}}\| V\|_r\geq 1.$
 More generally, let $u\in
W^{1,p}_0(D)$ be a weak solution of
\begin{equation}
-\Delta_pu=V(x)f(x,u,\nabla u).\label{pde}
\end{equation}
There are numerous works giving {\it
sufficient} conditions for the existence of solutions of
 equations of this form  with $D$ being the unit
interval, see \cite{OR}, \cite{BGG} and the references therein. In
\cite{BGG}, the authors prove the existence of multiple solutions
for a family of boundary value problems that include \rq{pde},
assuming that $f$ is continuous and non-negative, and $V$ is
continuous on $(0,1)$, does not vanish on any open subinterval, and
is $L^1$. The following result partly generalizes Theorem
\ref{Sbeta} and Theorem 3.7 in \cite{DEHL}.
\begin{Thm}\label{gennonlinS}
Assume that \rq{pde} holds, and
\begin{equation} |f(x,y,z)|\leq  |y|^{\beta +1}|z|^{\gamma},\label{abS}
\end{equation}
with constants $\beta \geq -1$ and $\gamma \geq 0$.
Assume $V\in L^r(D)$ with $p\neq n$ and
\begin{equation}
\frac{1}{r}+\frac{\beta +2}{\ov q}+\frac{\gamma }{p}= 1,\label{rbqpS}
\end{equation}
Then
\begin{equation}
 K_{p,\ov{q}}^{p-\gamma} \| V\|_r\|
u\|_{\ov{q}}^{2+\beta -p+\gamma}\geq 1.\label{msrnlS}
\end{equation}
This result also holds when $p=n$ and $\bar q$ is replaced by
$q<\infty$ in \rq{rbqpS} and \rq{msrnlS}.
\end{Thm}

{\it Proof: }  We begin with $p\neq n$. We will assume $\beta >-1$
and $\gamma > 0$. The proof for $\beta = -1$ is similar. The proof for
$\gamma =0$ is similar to the proof of \rq{betaS}.
By \rq{abS}, \rq{rbqpS} and $u\in
L^{\ov q}(D)$, we have $Vf(x,u,\nabla u)\in L^{(\ov{q})^*}$. We can
apply Lemma~\ref{GreenGG} and Holder's inequality to get
\begin{eqnarray}
\| \nabla u\|_p^p
& = & <Vf,u>\nonumber \\
&\leq & \| V\|_r \| f\|_t\| u\|_{\ov{q}},\label{part1}
\end{eqnarray}
where $t$ is defined by
$\frac{1}{r}+\frac{1}{t}+\frac{1}{\ov{q}}=1$. Define $j$ by $j(\beta
+1) t=\ov{q}$ for $p<n$ and $j=\infty$ for $p>n$. Let $k=p/\gamma
t$. Note that $\frac{1}{j}+\frac{1}{k}= t((\beta +1)/\ov{q}
+\gamma/p) = 1$, by \rq{rbqpS}. By Holder again and \rq{abS}, we get
\begin{eqnarray}
\| f\|_t& \leq & \| |u|^{\beta +1}|\nabla u|^{\gamma}\|_t \nonumber \\
&\leq & \| |u|^{(\beta +1)t}\|_j^{1/t}\| |\nabla u|^{\gamma
t}\|_k^{1/t}\nonumber \\
& =  &\| u\|_{\ov{q}}^{\beta +1}\| \nabla
u\|_{p}^{\gamma}.\label{part2}
\end{eqnarray}
Combining Sobolev's inequality, \rq{part1} and \rq{part2}, we get
$$\| u\|_{\ov{q}}^{p-\gamma} \le K_{p,\ov{q}}^{p-\gamma} \| \nabla u\|_p^{p-\gamma}
\leq  K_{p,\ov{q}}^{p-\gamma} \| V\|_r\| u\|_{\ov{q}}^{2+\beta}.\ \
$$
This proves \rq{msrnlS}.

For the case $p=n$, we cannot assume $u\in L^{\bar{q}}$, but we have $u\in L^q$ for all $q<\infty$. Assuming
\rq{rbqpS} holds with some finite $q$ replacing $\bar{q}$, the proof for this case is the same.
 $\Box$

\section*{ 5. Appendix}
\setcounter{section}{5} \setcounter{Thm}{0}
\setcounter{equation}{0}\setcounter{subsection}{0}

\subsection{Extremals and their Euler Lagrange equations}

\begin{Lemma}\label{SoboExtr2}
Let $p>n$ and $1\leq q\leq \infty$. Then there is a continuous
non-negative Sobolev extremal $u_* \in W^{1,p}_0 (D)$, with
$\|\nabla u_*\| _{p} = 1$ and
\begin{equation} \|u_*\|_{q}=K_{q,p}.
\label{extpn}
\end{equation}
\end{Lemma}

\noindent {\it Proof of Lemma \ref{SoboExtr2}:} We prove the result
for $q=\infty$. The proof for $q<\infty$ then follows by using
elementary arguments and observing that compact subsets of
$L^{\infty}(D)$ are compact subsets of $L^{q}(D)$. Let $B_W =
\{u\in W^{1,p}_0(D): \|\nabla u \|_{p} \leq 1\}$. For $p>1$,
$W^{1,p}(D)$ is reflexive. Since $W_0^{1,p}(D)$ is a closed subspace
of $W^{1,p}(D)$, it is also reflexive. Thus $B_W$ is weakly compact
with respect to the Sobolev norm. Moreover, the inclusion
$W^{1,p}_{0}(D) \to C^0(D)$ is compact. Let $\{u_n\}$ be a sequence
in $W_0^{1,p}(D)$ such that
 $$\lim_{n \rightarrow\infty}\frac{\|u_n\|_{\infty} }{\|  \nabla
 u_n\|_{p}} =K_{\infty ,p}.$$
We can assume by scaling that $\| \nabla u_n\|_{p} =1$. Since $B_W$
is weakly compact in  $W_0^{1,p}(D)$, there exists a subsequence $\{
u_{n_k}\} \subset \{u_n\}$  that  converges weakly to some $u_*\in
B_W$. By the compactness of the inclusion $W^{1,p}_0(D)\ra C^0(D)$,
there is a subsequence of $\{u_{n_k}\}$, that we label again with
 $\{ u_{n_k}\} $,  that converges  to some   $w\in C^0(D)$ in the strong topology of $C^0(D)$.  That is,
 $\dsize\lim_{k\to\infty} \| u_{n_k}- w\|_{\infty}=0.$
  But $u_{n_k}\to  u_*$ also in the weak topology of $ C^0(D)$, ie. pointwise
  and so $u_*=w$ a.e.;
 consequently, $w\in B_W$ and $\| \nabla w\|_{p} \leq 1$. We have
   $$
    K_{\infty ,p} =\lim_{k\to\infty} {\|u_{n_k}\|_{\infty}}=\|w\|_{\infty},
 $$
 and so
$\dsize \frac{\|w\|_{\infty} }{\|  \nabla
 w\|_{p}}
 \geq K_{\infty ,p}$. But recall that $w=u_*\in W^{1,p}_0(D)$, and so
$\dsize \frac{\|w\|_{\infty} }{\|  \nabla
 w\|_{p}}
 \leq K_{\infty ,p}$ proving \rq{extpn} for $u_*=w$. If $u_*$ is not already non-negative, we can
replace it by $|u_*|$,  with no effect on \rq{extpn} (see
\cite{LL}). $\Box$

\bigskip

Recall that $\bar{q}=\frac{np}{n-p}$ for $1<p<n$ and
$\bar{q}=\infty$ for $n\leq p$.
\begin{Lemma}\label{SoboExtr} Suppose $1<p\leq n$ and $1<q<\bar {q}$.
Then there is a non-negative Sobolev extremal
$u_* \in W^{1,p}_0 (D)$ with $\|\nabla u_*\| _{p} = 1$ and
$$\|u_*\|_q=K_{q,p}.$$
\end{Lemma}
The proof of this result is a straightforward adaptation of the proof of Lemma
 \ref{SoboExtr2} (also see \cite{DEHL}), and is left to the reader.

\medskip

\begin{Lemma}\label{EuLap} Let $1<p<\infty$ and
$p\leq q < \bar q$, and let $u_*$ be a Sobolev extremal as in Lemma
\ref{SoboExtr2} or Lemma \ref{SoboExtr}. Then \BEL{PDE} -\Delta_p
u_*= \frac{u_*^{q-1}}{ \| u_*\|_{q}^{q}}. \EE
\end{Lemma}

The proof of Lemma~\ref{EuLap} is similar to (\cite{DEHL}, Lemma 5.3), and is left to
the reader.

\subsection{Green's identities for divergence and Orlicz forms.}

\begin{Lemma}\label{GreenG}  Let $u\in W^{1,p}_0(D)$, with  $-{\rm div}(\nabla  u)= V u$ in the distribution sense. Assume either
\medskip
A) $p<n$ and $V \in L^{n/p} (D)$, or

B) $p\ge n$ and $V \in L^{r} (D)$ for some $r>1$, or

C) $n=p$, $u$ is real-valued and   $V\in L{\rm log}^{n-1}L(D)$. Then
\begin{equation}\label{greenh2}
\int_D  |\nabla  u|^p dx= \int_D V |u|^p dx.
\end{equation}
\end{Lemma}

\medskip
\noindent {\it Proof.} By assumption \BEL{pde4}   \int_D   |\nabla
u|^{p-2}\nabla u \cdot \nabla \psi    dx= \int_D V |u|^{p-2}u \psi \
dx \EE for every $\psi \in C^\infty_0(D)$.   Let $\{\psi_n\}$ be a
sequence of functions in $C^\infty_0(D)$ that converges to $\ov u$
in $W^{1,p}_0(D)$. Then
$$
\int_D|\nabla u|^{p-2}\nabla  u\cdot  \nabla \psi_n \,  dx -\int_D|\nabla  u|^p\, dx \leq ||  \nabla \psi_n - \nabla \ov u||_{p} |||\nabla u|^{p-1}||_{p^*} \to 0.
$$
To complete  the proof of \rq{greenh2}, it suffices to show that
$V|u|^{p-2}u\psi_n $ converges  to $V|u|^p$ in $L^{1}(D)$ in each case.
Assume A), that $p<n$ and $V \in L^{n/p} (D)$. By
Sobolev's inequality,  $\psi_n$ converges to $\ov{u}$ in
$L^{\frac{pn}{n-p}} (D)=L^{\ov{q}}(D)$. By Holder's inequality,
  $$\dsize \int_D |V|u|^{p-2} u \psi_n  -  V |u|^p|\, dx
  \leq   ||V  ||_{n/p}||u ^{p-1}||_{\ov{q}/(p-1)}
 \| \psi_n- \ov u\|_{\ov q}\to 0. $$

The proof in case B) is similar. So, assume C), and without loss of generality, that $\|\nabla
u\|_n= 1$. Let  $\{\psi_m\}\subset C_0^{\infty}(D)$ converge to $u$
in  $W^{1,n}_0(D)$. We can choose $\lambda_m \downarrow 0 $ so that
$\| \nabla (u-\psi_m) \|_n \lambda_m^{-1} \rightarrow 0$. The
Moser-Trudinger inequality implies, for $n$ large,
$$ \int_D
e^{\left(\alpha_n\frac{|u-\psi_m|}{\lambda_n}\right)^{n/(n-1)}}
dx\leq \int_D e^{\left( \alpha_n\frac{|u-\psi_m|}{\|\nabla (u-
\psi_m)\|_n}\right)^{n/(n-1)}}dx< (C_n+1)|D|<\infty,$$  with  $C_n$
independent of $u$ and $\psi_m$. A similar inequality holds when
$\frac{u-\psi_m}{\lambda_m}$ is replaced by $u$.

Let $M$ and $N$ be the functions defined in \rq{Mn}. We have $M(t)\leq Ce^{\alpha_n^{n/(n-1)} t^{1/(n-1)}}$ and $N(t)\sim t\log^{n-1}(t)$ for large $t$.  It follows from $V\in
L\log^{n-1}L(D)$ that
\begin{equation}
\int_DN(|V|)<\infty .\label{NO}
\end{equation}

 Using Young's inequality, the inequality $|ab|\leq
\frac{n-1}{n}a^{n/(n-1)}+\frac{1}{n}b^{n}$, and H\"older's inequality:
\begin{eqnarray*}
\int_D \left|u^{n-1} \,V (u-\psi_m )\right| dx & = & \lambda_m
\int_D |\frac{u^{n-1}
(u-\psi_m )}{\lambda_m}V| dx \\
& \leq &  \lambda_m \int_D M\left(\frac{|u|^{n-1}
(u-\psi_m )|}{\lambda_m}\right) dx+ \lambda_m\int_D N(|V|) dx \\
 & \leq & C\lambda_m \int_D \exp \left(\alpha_n^{\frac{n}{n-1}}|u|\left(\frac{
|u-\psi_m |}{\lambda_m}\right )^{1/(n-1)}\right )dx + \lambda_m\int_D N(|V|) dx\\
 & \leq & C\lambda_m \int_D
 e^{\frac{n-1}{n}|\alpha_nu|^{\frac{n}{n-1}}}
 e^{\frac{1}{n}(\alpha_n|u-\psi_m|\lambda_m^{-1})^{\frac{n}{n-1}}} dx+ \lambda_m\int_D N(|V|) dx\\
 & \leq & \lambda_m\left\{C\left(\int_D e^{|\alpha_nu|^{{\frac{n}{n-1}}}}dx\right)^{\frac{n-1}{n}}
 \left(\int_D e^{(\alpha_n\frac{|u-\psi_m|}{\lambda_m})^{{\frac{n}{n-1}}}}
 dx\right)^{\frac 1 n}  +\int_D N(|V|)dx\right\}
\end{eqnarray*}
Thus $ u \psi_m|u|^{n-2}V$ converges in $L^1(D)$ to $|u|^nV. \Box$

\bigskip

The next lemma is used in Section 4. It contains Lemma \ref{GreenG}
parts A) and B) as special cases.
\begin{Lemma}\label{GreenGG}   Let $u\in W^{1,p}_0(D)$ be a  solution in the
distribution sense of $-\Delta_pu= F$. Assume $F \in L^{(\ov q)^*}
(D)$ when $n\neq p$, and $F \in L^{r} (D)$ for some $r>1$ when $n=p$. Then,
\begin{equation}\label{gr3}
\int_D  |\nabla  u|^p dx= \int_D F\ov u dx.
\end{equation}
\end{Lemma}

\medskip
\noindent {\it Proof.}  We have $\int_D |\nabla u|^{p-2}\nabla u
\cdot \nabla \psi dx= \int_D F  \psi \ dx$ for every $\psi \in
C^\infty_0(D)$. The rest is similar to the proof of Lemma
\ref{GreenG}, part A).


\bigskip



 \begin{thebibliography}{999999}



\bibitem[AS]{AS} Astarita, G.; Marrucci, G., {\it Principles of Non-Newtonian Fluid
Mechanics}, McGraw-Hill, New York, (1974).

\bibitem[B]{B} Belaud, Y.; {\it Time-vanishing properties of solutions of some degenerate parabolic
 equations with strong absorption.}, Adv. Nonlinear Stud. 1, (2001), no. 2, p.117–152.

\bibitem[BGG] {BGG} Bai,Z., Gui, Z., and Ge, W. {\it Multiple positive
solutions for some $p$-Laplacian boundary value problems}, Math.
Anal. Appl. 300, (2004), p.477-490.



 \bibitem[CEP]{CEP} Cuccu, F., Emamizadeh, B., and Porru. {\it
 Optimization of the first eigenvalue in problems involving the
 $p$-Laplacian}, Proc. AMS 137, (2009), p.1677-1687.

\bibitem[DEHL]{DEHL}    De Carli, L; Edward, J; Hudson, S; Leckband, M. {\it  Minimal support results for  Schr\"odinger  equations},
  to appear in Forum Mathematicum.

\bibitem[DH1]{DH1}    De Carli, L; Hudson, S. {\it  Geometric Remarks on  the Level Curves of Harmonic Functions},
  Bull. London Math. Soc. 42  no. 1  (2010),  83--95.

\bibitem [DH2] {DH2}   De Carli, L.; Hudson, S. {\it A Faber-Krahn inequality for solutions of  Schr\"odinger's equations}, to appear in
Adv. Math.

\bibitem[GT]{GT} Gilbarg, D.; Trudinger, N. {\it Elliptic partial
differential equations of second order, 2nd edition},
Springer-Verlag, Berlin (1983).


\bibitem [H] {H} Hudson S.    {\it Remarks on extremals of minimal support inequalities}, to appear in "{\it Recent Advances in Harmonic Analysis and Applications (In Honor of Konstantin Oskolkov)}, Springer Proceedings in Mathematics (2012)



\bibitem [Ho]{Ho} Holopainen, I. {\it Quasiregular mappings and the p -Laplace operator. Heat kernels and analysis on manifolds, graphs, and metric spaces}, Contemp. Math. 338, Amer. Math. Soc.(2003), 219–239


\bibitem [HR] {HR} Rao R. and He, Q.
{\it Non-zero Solution for the Quasi-linear Elliptic Equation},
Chin. Quart. J. of Math. 24 no.1 (2009),  117-124.


\bibitem [KR] {KR} Krasnosel'ski\'{\i}, M. and Rutickii, Y.
{\it Convex Functions and Orlicz Spaces}, P. Noordhoff Ltd., Groningen (1961).


\bibitem [Le]{Le} L$\hat{e}$, A. {\it Eigenvalue problems for the $p$-Laplacian}, Non-linear Anal. 64 (2006), 1057-1099.









\bibitem [Li]{Li} Li, Yu Xiang.  {\it Remarks on the extremal functions for the Moser-Trudinger inequality} Acta Math. Sin. (Engl. Ser.) 22 (2006), no. 2, 545-550


\bibitem[LL] {LL}   Lieb, E. H.;  Loss, M. {\it Analysis}, 2nd edition, American Mathematical Society (2001).



\bibitem [M] {M} Moser, J.  {\it A Sharp form of an Inequality by N. Trudinger},
Indiana Univ. Math. 20  (1971), 1077--1092.


\bibitem [OR] {OR} O'Regan, D. {\it Some general existence principles and results for $ (f(y ' )) ' =qf(t,y,y ' ),0<t<1$} , SIAM J. Math. Anal. 24 (1993), no. 3, 648–668.



\bibitem[Ru]{Ru} Rudin, W.  {\it Real and complex analysis}, Third edition, McGraw-Hill   (1987).



\bibitem[St]{St} Struwe, M.  {\it Variational Methods: Applications to Nonlinear Partial Differential Equations and Hamiltonian Systems}, (4th edition), Springer(2008).

\bibitem [T] {T} Talenti, G.  {\it Best   Constant  in  Sobolev  Inequality}, Ann. Mat. Pura Appl.  110, no. 4  (1976), 353--–372.

\bibitem [U] {U} Uhlenbeck, K. {\it Regularity for a class of non-linear elliptic systems}. Acta. Math. 138 (1977) 219-240. (geometric application of p-lap)



\end{thebibliography}
\end{document}